\documentclass[12pt,draft]{amsart}
\usepackage{amsmath,amsthm,latexsym,amscd,amsbsy,amssymb}
 \relax

\chardef\bslash=`\\ 

\makeatletter
\def\verbatim{\interlinepenalty\@M \@verbatim
  \leftskip\@totalleftmargin\advance\leftskip2pc
  \frenchspacing\@vobeyspaces \@xverbatim}
\makeatother
\makeatletter
  \def\dgt@k{\dg@DX=-3 \dg@DY=2 \dg@SIZE=3} 
\makeatother

\makeatletter
  \def\dgt@kk{\dg@DX=3 \dg@DY=-1 \dg@SIZE=3}%
\makeatother

\theoremstyle{plain}
\newtheorem{thm}{Theorem}[section]

\newtheorem{lem}[thm]{Lemma}

\theoremstyle{definition}

\numberwithin{equation}{section}

\font\f=msbm10

\begin{document}


\title[Nonmetrizable $\text{ANR}$'s admitting a group structure are manifolds]{Nonmetrizable $\text{ANR}$'s admitting a group structure are manifolds}
\author{Alex Chigogidze}
\address{Department of Mathematical Sciences
University of North Carolina at Greensboro,
383 Bryan Building, Greensboro, NC, 27402,
USA}
\email{chigogidze@uncg.edu}



\begin{abstract}{It is shown that a nonmetrizable $\text{ANR}$-space of weight $\tau > \omega$, admitting a group structure, is (topologically) an $\text{\f R}^{\tau}$-manifold.}
\end{abstract}

\maketitle
\markboth{A.~Chigogidze}{Nonmetrizable $\text{ANR}$'s admitting a group structure are manifolds}

It is well known \cite[Corollary 1]{DT} that if a separable complete $\text{ANR}$ carries a topological group structure then either this is a Lie group or the $\text{ANR}$ is an $l_{2}$-manifold. We show that in the non-metrizable case the situation, is certain sense, is simpler. Namely we prove the following result.

\begin{thm}\label{T:main}
Let $\tau > \omega$ and $G$ be an $\text{ANR}$-space of weight $\tau$ admitting a group structure. Then $G$ is an $\mathbb{R}^{\tau}$-manifold.
\end{thm}

We refer the reader to \cite{chibook96} for comprehensive discussion of the general (non-metrizable) theory of absolute retracts, soft maps and spectral techniques. Here all we need to emphasize is that $\text{ANR}$-spaces are defined as retracts of functionally open subspaces of powers of the real line and that a map $p \colon X \to Y$ is soft iff there is a retraction $r \colon Y \times P \to X$, where $P$ is an $\text{AR}$-space, such that $pr = \pi_{Y}$. Two main properties of soft maps used below are: (a) soft maps are retractions (and hence admit sections); (b) inverse images of singletons are $\text{AR}$'s.

We begin with the following two statements needed in the proof of our result.

\begin{lem}\label{L:trivial}
Let $p \colon X \to Y$ be a continuous surjective homomorphism of topological groups. If there exists a continuous map $q \colon Y \to X$ such that $p q = \operatorname{id}_{Y}$, then $p$ is a trivial fibration with fiber $\operatorname{ker}(p)$. More precisely there exists a homeomorphism $h \colon X \to Y \times \operatorname{ker}(p)$ such that $\pi_{Y}h = p$.
\end{lem}
\begin{proof}
The required homeomorphism $h$ and its inverse $t \colon Y \times \operatorname{ker}(p) \to X$ are defined as follows:
\[ h(x) = \left( p(x), x\cdot q(p(x))^{-1}\right) , x \in X ;\;\;\; t(y,a) = a\cdot q(y), (y,a ) \in Y \times \operatorname{ker}(p).\]

First, note that $h$ is well defined. Indeed, for $x \in X$ we have
\begin{multline*}
 p\left( x\cdot q(p(x))^{-1}\right) = p(x)\cdot p\left( q(p(x))^{-1}\right) = p(x)\cdot \left( p(q(p(x))\right)^{-1} = \\ p(x) \cdot \left(p(x)\right)^{-1} = e_{Y} ,
\end{multline*}
which shows that $h(x) \in Y \times \operatorname{ker}(p)$.

For each $(y,a) \in Y \times \operatorname{ker}(p)$ we have

\begin{multline*}
h(a\cdot q(y)) = (p(a\cdot q(y)), (a\cdot q(y))\cdot (q(p(a\cdot q(y)))^{-1}) =\\ (p(a)\cdot p(q(y)), a\cdot q(y)\cdot (q(p(a) \cdot p(q(y)))^{-1}) = \left( e_{Y} \cdot y, a\cdot q(y)\cdot (q(e_{Y}\cdot y))^{-1}\right) =\\
 (y,a).
\end{multline*}

Also for each $x \in X$ we have

\[ t(h(x)) = t\left( p(x),x\cdot (q(p(x))^{-1}\right) = x\cdot (q(p(x))^{-1}\cdot q(p(x)) = x\cdot e_{X} = x .\] 

In other words $ht = \operatorname{id}_{Y \times \operatorname{ker}(p)}$ and $th = \operatorname{id}_{X}$. Thus $h$ is a homeomorphism.
\end{proof}

\begin{lem}\label{L:emb}
Let $\tau \geq \omega$. Every $\text{ANR}$-group of
weight $\tau \geq \omega$
is topologically and algebraically isomorphic to a closed
and $\text{C}$-embedded
subgroup of the product 
$\displaystyle \prod\{ G_{t} \colon t \in T\}$, where $|T| = \tau$ and each $G_{t}$,
$t \in T$, is a Polish $\text{AR}$-group homeomorphic to $\mathbb{R}^{\omega}$.
\end{lem}
\begin{proof}
This follows from \cite[Propositions 6.1.4, 6.5]{chibook96} and \cite[Proposition 4.1]{fund}.
\end{proof}

Now we prove our main result.

\bigskip

\noindent {\it Proof of Theorem \ref{T:main}}.
First we construct a well-ordered inverse spectrum
${\mathcal S}_{G} = \{ G_{\alpha},
p_{\alpha}^{\alpha +1},\tau \}$ satisfying the following properties:

\begin{enumerate}
\item
$G$ is topologically and algebraically isomorphic to $\lim {\mathcal S}_{G}$.
\item
for each $\alpha < \tau$, $G_{\alpha}$ is a $\text{ANR}$-group and $p^{\alpha +1}_{\alpha}
\colon G_{\alpha +1} \to G_{\alpha}$ is a soft
homomorphism such that $\operatorname{ker}\left( p_{\alpha}^{\alpha +1}\right)$ is a non-compact Polish $\text{AR}$-group.
\item
If $\beta < \tau$ is a limit ordinal, then the diagonal product 
\[\triangle\{ p_{\alpha}^{\beta} \colon \alpha < \beta\}
\colon G_{\beta} \to \lim\{ G_{\alpha},
p_{\alpha}^{\alpha+1}, \alpha < \beta \}\]
is a topological and algebraic isomorphism.
\item
$G_{0}$ is a Polish $\text{ANR}$-group. 
\end{enumerate}

By Lemma \ref{L:emb}, we may
assume that $G$ is a closed and $\text{C}$-embedded subgroup
of the product $\displaystyle X = \prod\{ X_{a} \colon a \in A\}$,
$|A| = \tau$, of Polish $\text{AR}$-groups $X_{a}$, $a \in A$. Since $G$ is an $\text{ANR}$-space, there exist a functionally open subspace $L$ of the product $X$ and a retraction $r \colon L \to G$. Choose a countable subset $A_{L} \subseteq A$ and an open subset $L_{G} \subseteq \prod\{ X_{a} \colon a \in A_{L}\}$ such that $L = L_{G} \times \prod\{ X_{a} \colon a \in A \setminus A_{L}\}$. 

Next let us denote by 
\[ \pi_{B} \colon \prod\{ X_{a} \colon a \in A\} \to
\prod\{ X_{a} \colon a \in B\}\] 
and
\[\pi^{B}_{C} \colon \prod\{ X_{a} \colon a \in B\} \to
\prod\{ X_{a} \colon a \in C \}\]
the natural projections onto the corresponding subproducts
($C \subseteq B \subseteq A$). Let also $G_{B}$ denote the subspace $\pi_{B}(G)$ of the product $X_{B} = \prod\{ X_{a} \colon a \in B\}$ and $L_{B} = L_{G} \times \prod\{ X_{a} \colon a \in B \setminus A_{L}\}$ for each $B \subseteq A$ with $A_{L} \subseteq B$.

We call a subset $B \subseteq A$ admissible if $A_{L} \subseteq B$ and the
following equality
\[ \pi_{B}(r(x)) = \pi_{B}(x)\]
is true for each point $x \in \pi_{B}^{-1}\left(G_{B}
\right)$. We need the following properties of admissible sets.

{\em Claim 1. The union of arbitrary collection of admissible
sets is admissible}. 

Indeed let $\{ B_{t} \colon t \in T\}$ be a collection of
admissible sets and $B = \cup \{ B_{t} \colon t \in T\}$.
Let $x \in \pi_{B}^{-1}\left(G_{B}\right)$. Clearly
$x \in \pi_{B_{t}}^{-1}\left(G_{B_{t}}\right)$ and we have 
\[ \pi_{B_{t}}(r(x)) = \pi_{B_{t}}(x) \;\;\text{for
each}\;\; t \in T .\]
Since $B = \cup \{ B_{t} \colon t \in T\}$ it follows that $\pi_{B}(x) = \pi_{B}(r(x))$.

{\em Claim 2. If $B$ is an admissible subset of $A$, then the
restriction $\pi_{B}|G \colon G \to G_{B}$ is a soft map (in the sense of \cite[Definition 6.1.12 ]{chibook96}). If $C$ and $B$ are admissible subsets with $C \subseteq B$, then $\pi_{C}^{B}|G_{B} \colon G_{B} \to G_{C}$ is also soft.}

Since the restriction of the projection $\pi_{B}$ onto $\pi_{B}^{-1}(G_{B})$ is a trivial fibration (with the fiber $\prod\{ X_{a} \colon a \in A \setminus B\}$) whose fiber is an $\text{AR}$-space (recall that each $X_{a}$ is an absolute retract and consequently so are their products), it follows that $\pi_{B}|\pi_{B}^{-1}(G_{B})$ is soft. The admissibility of $B$ implies that $\pi_{B}r|\pi_{B}^{-1}(G_{B}) = \pi_{B}|\pi^{-1}(G_{B})$. Since $r$ is a retraction it follows that $\pi_{B}|G$ is also soft (as a retract of $\pi_{B}|\pi_{B}^{-1}(G_{B})$). The second part of the claim follows from \cite[Lemma 6.1.15]{chibook96}.

{\em Claim 3. Each countable subset of $A$ is contained in a countable admissible subset of $A$}.

The inverse spectrum ${\mathcal S}(G)= \{ \operatorname{cl}_{L_{B}}\left( G_{B}\right), \pi_{C}^{B}|\operatorname{cl}_{L_{B}}\left( G_{B}\right), C,B \in \exp_{\omega}(A,A_{L})\}$ is
a factorizing (since $G$ is $C$-embedded in $X$) $\omega$-spectrum in the sense of \cite[Section 1.3.2]{chibook96} and $\lim {\mathcal S}(G) = G$. Similarly ${\mathcal S}(L) = \{ L_{B}, \pi_{C}^{B}|L_{B}, C, B \in \exp_{\omega}(A,A_{L})\}$ is also a factorizing $\omega$-spectrum and $\lim {\mathcal S}(L) = L$. Applying \cite[Theorem 1.3.6]{chibook96} to the map $r \colon L \to G$ we conclude that the collection $\exp_{\omega}(A, A_{L}$) contains a cofinal (and even $\omega$-closed) subcollection  ${\mathcal A}_{r}$ such that for each $B \in {\mathcal A}_{r}$ there exists a map $r_{B} \colon L_{B} \to \operatorname{cl}_{L_{B}}(G_{B})$ such that $\pi_{B}r = r_{B}\pi_{B}|L$. Note that for every such $B$, $G_{B}$ is closed in $L_{B}$ and $r_{B} \colon L_{B} \to G_{B}$ is a retraction. If $ x \in \pi_{B}^{-1}(G_{B})$ and $B \in {\mathcal A}_{r}$, then
\[ \pi_{B}(x) = r_{B}(\pi_{B}(x)) = \pi_{B}(r(x)) \]
which shows that every $B \in {\mathcal A}_{r}$ is admissible.

Since $|A| = \tau$ and $|A_{L}| = \omega$ we can
write $A \setminus A_{L} = \{ a_{\alpha} \colon \alpha < \tau\}$. By Claim 3,
each $a_{\alpha} \in A$ is contained in a countable admissible
subset $B_{\alpha} \subseteq A$. Let
$A_{\alpha} = \cup\{ B_{\beta} \colon \beta \leq \alpha\}$.
We use these sets to define a transfinite inverse spectrum
${\mathcal S} = \{ G_{\alpha}, p_{\alpha}^{\alpha +1}, \tau \}$
as follows. 

Let $A_{0}$ be a countable admissible subset of $A$ containing $A_{L} \cup \{ a_{0}\}$ and $G_{0} = G_{A_{0}}$.

Assume that the admissible subsets $A_{\beta} \subseteq A$ have already been constructed for each $\beta <\alpha$, $\alpha < \tau$, so that $a_{\beta} \in A_{\beta}$ for $\beta < \alpha$, $A_{\gamma} \subseteq A_{\beta}$ for $\gamma <\beta$ and the spaces $G_{\beta} = G_{A_{\beta}}$ and the homomorphisms $p_{\gamma}^{\beta} = \pi_{A_{\gamma}}^{A_{\beta}}|G_{\beta} \colon G_{\beta} \to G_{\gamma}$ satisfy the required properties (2) and (3).

If $\alpha$ is a limit ordinal, then let $A_{\alpha} = \cup\{ A_{\beta} \colon \beta < \alpha \}$ and $G_{\alpha} = G_{A_{\alpha}}$. 

If $\alpha = \beta +1$ then proceed as follows. Consider the fiber $\operatorname{ker}\left( \pi_{\beta}|G\right)$ of the homomorphism $\pi_{\beta}|G \colon G \to G_{\beta}$. This fiber is an $\text{AR}$-space (as an inverse image of a point under a soft map; see Claim 2). Since it is also a topological group we conclude that $\operatorname{ker}\left( \pi_{\beta}|G\right)$ is non-compact (much more general statement is true -- every compact $\text{AE}(4)$-group is trivial \cite{chi972}) -- otherwise it must be a singleton in which case the homomorphism $\pi_{\beta} \colon G \to G_{\beta}$ must be a homeomorphism. But this is impossible since $w(G) > w(G_{\beta})$. Consequently there exists a countable subset $C \subseteq A$ such that $\pi_{C}\left(\operatorname{ker}\left( \pi_{\beta}|G\right)\right)$ is also non-compact. By claim 3, there exists a countable admissible subset $B$ such that $C \cup \{ a_{\alpha}\} \subset B$. We now let $A_{\alpha} = A_{\beta} \cup B$, $G_{\alpha} = G_{A_{\alpha}}$ and $p_{\beta}^{\beta +1} = \pi_{A_{\beta}}^{A_{\alpha}}|G_{A_{\alpha}}$.
Note that $\operatorname{ker}\left( p_{\beta}^{\beta +1}\right)$ is a non-compact Polish $\text{AR}$-group. 

This completes the inductive step and the construction of the spectrum ${\mathcal S}_{G}$. Conditions (2) and (3) are satisfied by construction. Condition (1) can be insured by choosing the sets $A_{\alpha}$ so that $A = \cup\{ A_{\alpha} \colon \alpha < \tau \}$ which in turn is possible since the collection of all countable admissible sets in cofinal in $\exp_{\omega}(A,A_{L})$ (see the inductive step above).

The straightforward transfinite induction coupled with Lemma \ref{L:trivial} shows that $G$ is homeomorphic to the product $G_{0} \times \prod\{ \operatorname{ker}\left( p_{\alpha}^{\alpha +1}\right) \colon \alpha < \tau \}$. Since $p_{\alpha}^{\alpha +1}$ is a non-compact Polish $\text{AR}$-space for each $\alpha < \tau$, it follows (by \cite[Theorem 5.1]{tor}) that all countable (but infinite) subproducts of the product $\prod\{ \operatorname{ker}\left( p_{\alpha}^{\alpha +1}\right) \colon \alpha < \tau \}$ are homeomorphic to $\mathbb{R}^{\omega}$. Thus, $\prod\{ \operatorname{ker}\left( p_{\alpha}^{\alpha +1}\right) \colon \alpha < \tau \} \approx \mathbb{R}^{\tau}$ and $G \approx G_{0} \times \mathbb{R}^{\tau}$. Since $G_{0}$ is a Polish $\text{ANR}$-space, it follows from \cite[Corollary 7.3.4]{chibook96} that $G$ is an $\mathbb{R}^{\tau}$-manifold.



\providecommand{\bysame}{\leavevmode\hbox to3em{\hrulefill}\thinspace}



\end{document}